\pdfoutput=1
\documentclass[final]{amsart}
\usepackage{amsfonts}
\usepackage{amsmath}
\usepackage{amssymb}
\usepackage{amsthm}
\usepackage{graphicx}
\usepackage{graphics}

\theoremstyle{plain}
\newtheorem{proposition}{Proposition}[section]
\newtheorem{theorem}[proposition]{Theorem}
\newtheorem{lemma}[proposition]{Lemma}
\newtheorem{corollary}[proposition]{Corollary}

\newtheorem{conjecture}[proposition]{Conjecture}

\newtheorem*{HFW}{Hoffman-Wielandt Theorem}
\newtheorem*{WST}{Whitney Synthesis Theorem}

\DeclareMathOperator{\Rank}{Rank}
\DeclareMathOperator{\Ker}{Ker}
\DeclareMathOperator{\Max}{max}
\DeclareMathOperator{\Min}{min}
\DeclareMathOperator{\Nullity}{Nullity}

\newcommand{\RR}{\mathbb{R}}

\begin{document}
\author{Matthew Jacobs}
\title{Connecting Global and Universal Rigidity}
\thanks{This paper was supported by NSF grant DMS-0739392}

\begin{abstract}
A $d$-dimensional framework is an embedding of the vertices and edges of a graph in $\RR^d$.  A $d$-dimensional framework is globally rigid if every other $d$-dimensional framework with the same edge lengths has the same pairwise distances between the vertices. A graph is generically globally rigid in dimension $d$ ($d$-GGR) if every generic framework is globally rigid. The $d$-dimensional framework of a $d$-GGR graph is universally rigid if for $d' \geq d$ every $d'$-dimensional framework with the same edge lengths has the same pairwise distances between the vertices. We establish a strong connection between global and universal rigidity by showing that all 1 and 2-GGR graphs and an infinite number of higher dimensional $d$-GGR graphs have a generic universally rigid framework.       
\end{abstract}

\maketitle

\section{Introduction}

An interesting and difficult problem in graph theory is the graph embedding problem. The problem is as follows: given a collection of vertices in $\RR^d$ and the pairwise distances between some of them, find the positions of the vertices up to some arbitrary rigid motion. The problem has shown to be NP-HARD \cite{sx}; however, solutions may be found using semidefinite programming \cite{ll}. Unfortunately the embeddings given by semidefinite programming may be in some larger dimensional space rather than $d$-dimensional.  However, if the framework is \emph{universally rigid} then solutions are automatically constrained to be $d$-dimensional. This paper seeks to address the question raised by Gortler and Thurston:  ``For a graph $G$ which is $d$-GGR is there always a generic framework in the space of $d$-dimensional frameworks, $C^d(G)$  that is universally rigid?''\cite{gt}. This question is answered in the affirmative for $\RR$ and $\RR^2$, and for an infinite number of $d$-GGR graphs when $d>2$. 

The key tool to obtain this result is the Hennenberg operation.  The Hennenberg operation is a vertex and edge addition to a graph, which preserves important rigidity properties.  Connelly proved that Hennenberg operations preserve global rigidity \cite{cn2}. This paper expands that result by proving that Hennenberg operations preserve universal rigidity for suitable frameworks. Hennenberg operations and edge additions can be used to build $d$-GGR graphs from the complete graph $K_{d+2}$. All graphs built in this manner will possess a generic $d$-dimensional universally rigid ($d$-GUR) framework. Every non-trivial 1 and 2-GGR graph is the result of Hennenberg operations and edge additions to $K_{d+2}$, allowing us to conclude the main result.  
\begin{theorem} All 1 and 2-GGR graphs and an infinite number of $d$-GGR graphs for $d>2$ have a $d$-GUR framework. \end{theorem}

The complete answer to Gortler and Thurston's question is still unknown in higher dimensions but we offer the following conjecture. 

\begin{conjecture} All d-GGR graphs have a d-GUR framework. \end{conjecture}

Thanks to Kiril Ratmanski and Andrew Fanoe for their helpful comments and ideas. Thanks to Timothy Sun for showing that 1-GGR graphs can be characterized by the Whitney Synthesis Theorem. Special thanks to Dylan Thurston for his enormous help throughout the entire process.

\section{Background}
\textbf{Graphs and frameworks:}  A \emph{graph} $G$ is a set $\mathcal{V}$ of $v$ vertices, and a set $\mathcal{E}$ of $e$ edges containing some collection of 2 member subsets of $\mathcal{V}$.   A \emph{framework} $G(p)$, is an embedding of the vertices and edges in $\RR^d$ $\emph{p}=(p_1,p_2,\ldots,p_v)$ where each $p_i$ is the coordinate of a vertex in $\RR^d$.  Let $C^d(G)$ denote the set of all frameworks of $G$ in dimension $d$. A framework is \emph{generic} if the coordinates of its vertices do not satisfy any nontrivial polynomial with rational coefficients.

\textbf{Local Rigidity:} Two different embeddings $G(p)$, and  $G(q)$ of a graph $G$ are said to be \emph{equivalent} written $G(p)$$\sim$$G(q)$ if $ \lVert p_i - p_j\rVert^2=\lVert q_i - q_j\rVert^2$ for all $\{i,j\}$ edges in $\mathcal{E}$.   A framework is said to be \emph{locally flexible} if there exists a non-constant analytic path $g: [0,1]\to \RR^{vd}$ such that $g(0)=G(p)$ and $g(t)\sim G(p)  \;  \forall t \in [0,1]$ and $g$ cannot be extended to a Euclidean motion of $\RR^d$.  If a framework is not locally flexible then it is \emph{locally rigid}. 

\textbf{Infinitesimal Rigidity:} The \emph{half edge length squared function} $f: C^d(G) \to \RR^e$, where $f(G(p)) = \frac{1}{2}(\ldots, \lVert p_i - p_j\rVert^2, \ldots)$  $\forall \{i,j\} \in \mathcal{E}$ is a map from frameworks to the space of edge lengths.  If two frameworks are equivalent $G(p)\sim G(q)$ then by definition $f(G(p))=f(G(q))$. The Jacobian $df$ of the edge-length squared function is called the \emph{rigidity matrix}, and gives a simple way to determine local rigidity for generic frameworks. The matrix is a $e \times vd$ matrix where each row represents an edge and each column represents the one of the $d$ coordinates of a vertex. The rigidity matrix has an important connection to local rigidity, as the first derivatives of paths preserving edge lengths are in the kernel of the matrix.  Thus every Euclidean motion is in the kernel, and as long as the framework does not lie in some lower dimensional affine subspace of $\RR^d$, the dimension of the kernel is at least $\binom{d+1}{2}$.  A framework is \emph{infinitesimally rigid} if the kernel of the rigidity matrix is minimal, or equivalently $\Rank(df(G(q)))=vd - \binom{d+1}{2}$. If a framework is infinitesimally rigid then it is locally rigid. 

\begin{theorem}[Asimow and Roth \cite{ar}]
If a generic framework of a graph with d+1 or more vertices is locally rigid in $\RR^d$, then it is infinitesimally rigid.  Furthermore if there exists a single infinitesimally rigid framework, then every generic framework is infinitesimally rigid. \end{theorem}

\textbf{Stresses:}  An \emph{equilibrium stress} on a framework is an assignment of a real number $\omega_{ij}$ to each edge of a framework satisfying $\sum_{j} \omega_{ij}(p_i - p_j)=\mathbf{0}$. The \emph{stress vector} $\boldsymbol{\omega}$ is a $e \times 1$ vector consisting of the ordered stresses on each edge. It is not difficult to see that $\boldsymbol{\omega}$ is in the kernel of the transpose of the rigidity matrix, and in fact the space of all stresses is precisely this kernel.  An \emph{equilibrium stress matrix} $\Omega$ is a $v \times v$ matrix satisfying: 

(1) $\Omega_{i,j}=\Omega_{j,i}$

(2) $\Omega_{i,j}=0$ if $i \ne j$ and $\{i,j\} \notin \mathcal{E}$

(3) $\Omega_{i,j}=-\omega_{ij}$ if $i \ne j$ and $\{i,j\} \in \mathcal{E}$

(4) $\Omega_{i,i}$ is chosen so the row sum is zero. 

The entries of the stress matrix are uniquely determined by $\boldsymbol{\omega}$ so this paper will freely interchange between equilibrium stresses and stress matrices.  Since the row sum of the matrix is zero the vector of all ones $\mathbf{1}$ is always in the kernel of the stress matrix.  Additionally the $d$ coordinate projections of the vertices lie in the kernel of the stress matrix. This condition is easy to check since $\sum_{j} \omega_{ij}(p_i - p_j)=\mathbf{0}$ must be satisfied. A generic framework with $d+1$ or more vertices cannot lie in a lower dimensional affine subspace, so the $d$ coordinate projections are linearly independent.  Thus for a generic framework with $d+1$ or more vertices, the nullity of the stress matrix is at least $d+1$. 

\textbf{Redundant Rigidity:} An edge in a locally rigid framework is \emph{redundant} if removal of said edge preserves local rigidity.  A framework is \emph{redundantly rigid} if every edge of the framework is redundant.  
         
\begin{proposition}[see e.g. Frank and Jiang \cite{fj}] A framework is redundantly rigid if and only if there exists a stress which is nonzero on every edge. \end{proposition}

\textbf{Global Rigidity:} Two frameworks $G(p)$ and $G(q)$ are \emph{congruent} written $p\equiv q$ if $ \lVert p_i - p_j\rVert^2=\lVert q_i - q_j\rVert^2$ for all $i,j \in \mathcal{V}$.  A $d$-dimensional framework $G(p)$ is \emph{globally rigid in dimension} $d$ if every framework in $C^d(G)$ equivalent to $G(p)$ is also congruent to $G(p)$. 

\begin{theorem}[Connelly \cite{cn2}; Gortler, Healy, Thurston \cite{ght}]
A graph is generically globally rigid in dimension $d$ ($d$-GGR) if and only if it has a generic framework with a stress matrix of minimal nullity d+1, or it is the complete graph on d+1 or fewer vertices.  \end{theorem}

Connelly showed that the stress condition was sufficient, and Gortler, Healy, and Thurston showed that it was necessary, implying that global rigidity is a generic property. Since global rigidity is a generic property, it can be thought of as a property of a graph, and such graphs are called \emph{generically globally rigid in dimension} $d$ ($d$-GGR).  The stress condition is not particularly intuitive, and other attempts have been made to characterize global rigidity combinatorially.  

\textbf{Vertex connectivity:} A graph is $n$\emph{-vertex-connected} if the deletion of any $n-1$ vertices leaves a connected graph.  

\begin{theorem}[Hendrickson \cite{he}]
If a graph is $d$-GGR then it is redundantly rigid and d+1-vertex-connected. \end{theorem} 

In the same paper Hendrickson conjectured that these conditions would also be sufficient. In $\RR$ 2-vertex-connectedness subsumes redundant rigidity, and it is easy to see that all 2-vertex-connected graphs are 1-GGR. In $\RR^2$, Connelly and Jackson and Jord\'an proved that these conditions are sufficient for global rigidity (see Theorem 5.7) \cite{cn2}, \cite{jj}.  Unfortunately, in $\RR^3$ and above Connelly showed that there are graphs which satisfy Hendrickson's conditions, but fail to be globally rigid. In fact there is not even a combinatorial characterization of local rigidity in dimensions 3 or higher \cite{cn}.  The lack of a general combinatorial characterization of $d$-GGR graphs for $d>2$ makes Conjecture 1.2 difficult to resolve.  
  
\textbf{Universal Rigidity:}  A $d$-dimensional framework, $G(p)$ is \emph{universally rigid} if every framework in $C^{d'}(\mathcal{V})$ equivalent to $G(p)$ is also congruent to $G(p)$ where $d'\geq d$.  If the framework is also generic then it is a \emph{$d$-dimensional generic universally rigid framework} ($d$-GUR framework).  Unlike global rigidity, universal rigidity is not a generic property.  A graph $G$ is \emph{$d$-semi-universally rigid} ($d$-SUR) if it has a $d$-GUR framework, but has other generic frameworks which fail to be universally rigid in $d$-dimensions.  Figure 1 has a very simple example of a graph with this behavior. In section 6 we construct a non-trivial family of $d$-SUR graphs.  Certain classes of graphs, like complete graphs, are generically universally rigid.

\begin{figure}[htb]
\begin{center}
\leavevmode
\includegraphics[width=0.8\textwidth]{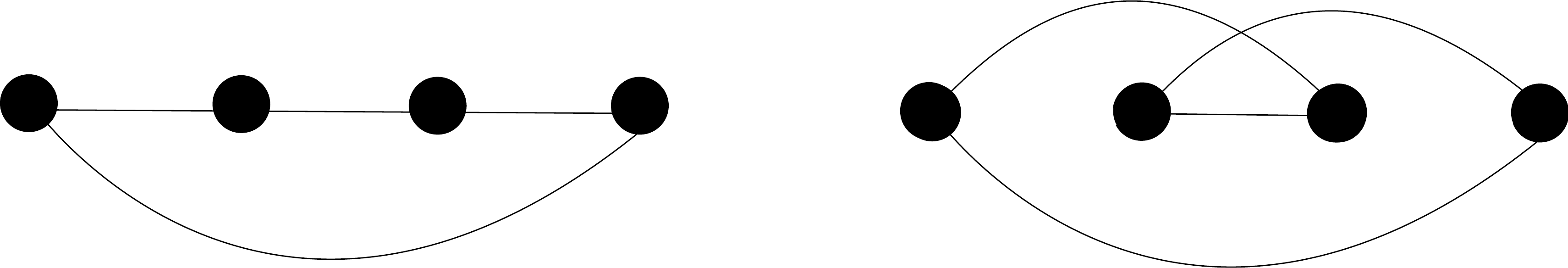}
\end{center}
\caption{Two different one dimensional frameworks of the same graph.  The framework on the left is universally rigid while the framework on the right is not.}
\label{figure1}
\end{figure}

\begin{proposition} Any framework of a complete graph is universally rigid. \end{proposition}
 
\begin{proof} For frameworks of complete graphs equivalence and congruence are the same condition. \end{proof}

\begin{proposition} If a framework is $d$-GUR then it is a framework of a $d$-GGR graph. \end{proposition}

\begin{proof} If a framework is universally rigid then it must be globally rigid in dimension~$d$.  Since the framework is generic all other generic $d$-dimensional frameworks are globally rigid so the graph is $d$-GGR.  \end{proof}

\textbf{Hennenberg Operations:} A $d$-Hennenberg operation is an edge and vertex addition to a graph that preserves important rigidity properties.  A $d$-Hennenberg operation is performed on a graph by deleting an edge $\{x,y\}$, adding a new vertex $z$ and edges $\{x,z\}$ and $\{y,z\}$, and adding an additional $d-1$ edges to $z$ which are connected to any $d-1$ distinct vertices that are not $x$ or $y$.  It is important to note that $d$-Hennenberg operations may only be performed on graphs with at least $d+1$ vertices, since the operation requires $d+1$ distinct vertices.  Connelly showed that $d$-Hennenberg operations on a $d$-GGR graph resulted in another $d$-GGR graph \cite{cn2}.  This paper extends that result to show that a $d$-Hennenberg operation takes a graph with a $d$-GUR framework to another graph with a $d$-GUR framework.

\section{Characterizing Universal Rigidity}

\begin{theorem}[Gortler and Thurston \cite{gt}] If a generic d-dimensional framework with d+2 or more vertices is universally rigid then it has a positive semi-definite (PSD) stress matrix of nullity d+1.  \end{theorem}
Alfakih proved the theorem for a matrix of any rank \cite{alf1}, \cite{alf2}. Gortler and Thurston strengthened the theorem by showing the stress matrix has minimal nullity. 

\textbf{Conics at infinity:} Let $Q$ be a symmetric $d \times d$ non-zero matrix.  A conic at infinity $C(Q)$ is the set of points $C(Q)=\{x \in \RR^d$ : $x^TQx=0\}$.  The edge directions of a framework $G(p)$ lie on a conic at infinity if there exists a non-zero symmetric matrix $Q$ satisfying $(p_i-p_j) \in C(Q)$ for all $\{i,j\} \in \mathcal{E}$.   

\begin{theorem}[Connelly] A d-dimensional framework with at least d+2 vertices is universally rigid if it has a PSD stress matrix of nullity d+1 and its edge directions do not lie on a conic at infinity. \end{theorem}
Connelly proved this result across a number of papers \cite{cn3} and \cite{cn4}, however due to the centrality of the result to this paper a condensed and slightly modified proof is presented below.  

\begin{proof} Let $E^{\boldsymbol{\omega^r}}_G$ be the energy function $E^{\boldsymbol{\omega^r}}_G: C^d(G) \to \RR$. \begin{equation} E^{\boldsymbol{\omega^r}}_G(p)=\sum_{\{i,j\}\in \mathcal{E}}\omega^{r}_{i,j}\lVert p_i - p_j\rVert^2, \end{equation} where $\omega^{r}_{ij}$ is the stress given by the fixed stress vector $\boldsymbol{\omega^r}$ for the framework $G(r)$.  Let $p^m$ be the $v \times 1$ vector of the $m^{th}$ coordinate entry of each of the $v$ vertices of $G(p)$, then an equivalent version of $E^{\boldsymbol{\omega^r}}_G$ is given by \begin{equation} E^{\boldsymbol{\omega^r}}_G(p)=\sum_{m=1}^d (p^m)^T\Omega_{G(r)} p^m, \end{equation} where again $\Omega_{G(r)}$ is the fixed stress matrix determined by $\boldsymbol{\omega^r}$.  Now let $G(p)$ be a framework in dimension $d$ with a PSD stress matrix $\Omega_{G(p)}$, and whose edges do not lie on a conic at infinity.  Let $G(q)$ be some framework in dimension $d'\geq d$, with $G(p)$$\sim$$G(q)$.  Now consider $E^{\boldsymbol{\omega^p}}_G$, where $\boldsymbol{\omega^p}$ is the stress vector determined by $\Omega_{G(p)}$. Since the coordinate projections of $G(p)$ are in the kernel of $\Omega_{G(p)}$, by equation (2) $E^{\boldsymbol{\omega^p}}_G(p)=0$.  By equation (1) $E^{\boldsymbol{\omega^p}}_G$ is a function on the edge lengths.  Thus if $G(p)$$\sim$$G(q)$ then $E_G(p)=E_G(q)=0$. $\Omega_{G(p)}$ is a PSD matrix, so $\exists C$ such that $C^TC=\Omega_{G(p)}$. Thus, \begin{equation} E_G(q)=\sum_{m=1}^{d'} (q^m)^TC^TCq^m=\sum_{m=1}^{d'} \lVert Cq^m\rVert^2=0 \end{equation} \begin{equation} Cq^m=0 \qquad \Rightarrow \qquad \Omega_{G(p)} q^m=0 \qquad \forall m: 1 \leq m \leq d. \end{equation} The kernel of $\Omega_{G(p)}$ consists of the coordinate projections $p^m$ and the vector of all ones $\mathbf{1}$, thus each coordinant projection $q^m$ must be some linear combination of these vectors. Therefore $q$ is obtained from $p$ by some affine transformation of $\RR^d$ to $\RR^{d'}$. In other words $q_i=Ap_i + \mathbf{b}$ for some fixed $d'\times d$ matrix $A$ and  a fixed $d' \times 1$ vector $\mathbf{b}$ for all verticies $i \in \mathcal{V}$. The transformation must preserve edge lengths so \begin{equation} (p_i-p_j)^T(p_i-p_j)=(Ap_i-Ap_j)^T(Ap_i-Ap_j)=(p_i-p_j)^TA^TA(p_i-p_j)\end{equation} \begin{equation} 0=(p_i-p_j)^T(I_{d}-A^TA)(p_i-p_j). \end{equation} $I_{d}-A^TA$ is clearly a symmetric matrix; thus either it is the zero matrix or the edge directions lie on a conic at infinity. Since the edge directions do not lie on a conic at infinity by hypothesis, $I_{d}=A^TA$.  Let $i$,$j$ be vertices not necessarily sharing an edge. Then \begin{equation} \lVert q_i-q_j\rVert^2=\lVert A(p_i-p_j)\rVert^2=(p_i-p_j)^TA^TA(p_i-p_j)=\lVert p_i-p_j\rVert^2. \end{equation}
Thus every framework $G(q)$ in dimension $d' \geq d$ equivalent to $G(p)$ is also congruent to $G(p)$ as desired. \end{proof}

\begin{proposition}[Connelly \cite{cn2}] Let G(p) be a generic d-dimensional framework of a finite graph G where each vertex has degree at least d. Then the edge directions of G(p) do not lie on a conic at infinity. \end{proposition} 

\begin{corollary} A generic d-dimensional framework with d+2 or more vertices is universally rigid if and only if it has a PSD stress matrix of nullity d+1. \end{corollary}
\begin{proof} Theorem 3.1 already establishes one direction.  To prove the other direction note the nullity of the stress matrix is $d+1$, therefore the framework is globally rigid and thus $(d+1)$-connected.  Since the graph is $(d+1)$-connected the degree of each vertex must be at least $d$.  By Proposition 3.3, the edge directions of the framework do not lie on a conic at infinity.  The framework then satisfies the hypothesis of Theorem 3.2. \end{proof}

\section{Matrix Perturbations and Kernels}

In this section we derive a few simple results about stress matrices that we will need for our main argument.  The first four results deal with perturbations of eigenvalues and the last two results relate to kernels.

\textbf{Normal matrices:} A real valued $n \times n$ matrix $A$ is \emph{normal} if $AA^T=A^TA$. It is easy to see that any symmetric matrix $S$ is normal since $S=S^T$.

\begin{HFW}[see e.g. Horn and Johnson, Section 63  \cite{hjm}] Let A, F be $n \times n$ matrices, with A, A+F normal.  Let $(\lambda_{1}, \lambda_{2}, \ldots, \lambda_{n})$ be the eigenvalues of A in some order and $(\mu_{1}, \mu_{2}, \ldots, \mu_{n})$ be the eigenvalues of A+F in some order.  Then there exists a permutation $\sigma(i)$ of the integers $1,2,...,n$ such that \begin{equation} \left[ \sum_{i=1}^n |\mu_{\sigma(i)}-\lambda_{i}|^2 \right]^{1/2} \! \! \! \leq \lVert F\rVert_{2} \end{equation} \end{HFW}

\begin{corollary} Let A and C be $n \times n$ symmetric matrices. Let B be the $n \times n$ matrix such that $A+B=C$. Let $|\lambda_{m}|$ be the smallest absolute value of an eigenvalue of A.  If $\lVert B \rVert_{2} < |\lambda_{m}|$ then C has at least as many positive and negative eigenvalues as A.  \end{corollary}

\begin{proof} $A$ and $C$ are symmetric, so they are normal and all their eigenvalues are real. Let $(\lambda_{1}, \lambda_{2}, \ldots, \lambda_{n})$ be the eigenvalues of $A$, and $(\mu_{1}, \mu_{2}, \ldots, \mu_{n})$ be the eigenvalues of $C$. By the Hoffman and Wielandt theorem, \begin{equation}\sum_{i=1}^n |\mu_{\sigma(i)}-\lambda_{i}|^2 < |\lambda_{m}|^2 \end{equation} \begin{equation} |\mu_{\sigma(i)}-\lambda_{i}|< |\lambda_{m}| \quad \forall i:  1 \leq i \leq n. \end{equation} Thus $\lambda_{i}-|\lambda_{m}| < \mu_{\sigma(i)} < |\lambda_{m}| + \lambda_{i}$.  Since $|\lambda_{i}| \geq |\lambda_{m}|$ if $\lambda_i>0$ then $\mu_{\sigma(i)}>0$ and if $\lambda_i<0$ then $\mu_{\sigma(i)}<0$.  Thus $C$ has at least as many positive and negative eigenvalues as $A$.  \end{proof}                 

\begin{corollary} Let A be a $v \times v$ stress matrix of nullity $d+1$ for some d-dimensional framework G(p). Let C be a $v \times v$ stress matrix for some d-dimensional framework H(q), not necessarily different from G(p). Let B be the matrix such that $A+B=C$. Let $|\lambda_{m}|$ be the smallest absolute value of an eigenvalue of A.  If $\lVert B \rVert_{2} < |\lambda_{m}|$ then A and C have the same number of positive and negative eigenvalues.  \end{corollary}

\begin{proof} $A$ and $C$ satisfy Corollary 4.1 so $C$ has at least as many positive and negative eigenvalues as $A$. Thus $C$ has nullity at most $d+1$.  But $d+1$ is the minimal nullity for stress matrices so the rest of the of the eigenvalues of $C$ must be zero.  Therefore $A$ and $C$ have the same number of positive and negative eigenvalues.  \end{proof}

\begin{corollary} Let A, B be $v \times v$ stress matrices for some d-dimensional framework G(p), where A is of nullity d+1.  Then $\exists \epsilon>0$ such that $A$ and $A+\epsilon B$ have the same number of positive and negative eigenvalues, and $A+\epsilon B$ is a stress matrix for G(p). \end{corollary} 

\begin{proof} It is easy to check that if $A$ and $B$ are stress matrices for a framework then $c_1A+c_2B, \; c_1, c_2 \in \RR$ is also a stress matrix for the same framework. Let $k=\lVert B \rVert_2$, and let $|\lambda_m|$ be the smallest absolute value of an eigenvalue of $A$.    Pick $\epsilon \leq \frac{|\lambda_m|}{2k}$ and apply Corollary 4.2. \end{proof}   

Now we turn to analyzing kernels of matrices. First we define a way to measure distances between linear subspaces. 

\textbf{Hausdorff Distance:} Let ($X$, $\rho$) be a metric space and let $A$, $B$ be subsets of $X$.  The Hausdorff distance between $A$ and $B$, written $\rho_H(A,B)$ is defined as \begin{equation} \rho_H(A,B)=\Max\{\inf_{a \in A} \sup_{b \in B} \rho(a,b), \, \inf_{b \in B} \sup_{a \in A} \rho(a,b)\}. \end{equation} The distance between linear subspaces $L_1$ and $L_2$ of $\RR^d$ is defined as the Hausdorff distance between $L_1 \cap \mathbb{S}^{d-1}$ and $L_2 \cap \mathbb{S}^{d-1}$, where $\mathbb{S}^{d-1}$ is the unit sphere in $\RR^{d}$. $\mathbb{S}^{d-1}$ is a bounded set so the distance will be well defined.

\begin{lemma}[Connelly \cite{cn2}] Suppose that $A(x_1, x_2, \ldots, x_n)$ is a matrix whose entries are integral polynomial functions of the real variables $x=(x_1,x_2, \ldots, x_n)$. Let $m$ be the maximum rank of $A(x)$ and suppose that $\Rank(A(x'))=m$.  Then $\forall \epsilon>0 \; \exists \delta>0$ such that for $|x-x'|<\delta$ the Hausdorff distance between the linear subspaces $\Ker(A(x))$ and $\Ker(A(x'))$ is less than $\epsilon$.   \end{lemma}

\begin{lemma} If A, B are PSD matrices then $\Ker(A + B)=(\Ker(A)\cap \Ker(B))$. \end{lemma}
\begin{proof} Let $M$, $N$ be matrices such that $A=M^TM$ and $B=N^TN$. Let $\boldsymbol{v} \in \Ker(A+B)$.  
\begin{equation} \boldsymbol{v}^T(A + B)\boldsymbol{v}=\boldsymbol{v}^TM^TM\boldsymbol{v} + \boldsymbol{v}^TN^TN\boldsymbol{v} = \lVert M\boldsymbol{v}\rVert^2 + \lVert N\boldsymbol{v}\rVert^2=0 \end{equation} 
Therefore, $M\boldsymbol{v}=0$ and $N\boldsymbol{v}=0$, which holds if and only if $A\boldsymbol{v}=0$ and $B\boldsymbol{v}=0$. Thus $\boldsymbol{v} \in \Ker(A+B)$ if and only if $\boldsymbol{v} \in \Ker(A)$ and $\boldsymbol{v} \in \Ker(B)$. It follows that $\Ker(A+B)=(\Ker(A) \cap \Ker(B))$. \end{proof}

\section{The Main Result}

\begin{proposition} Let H be obtained from G by a single edge addition.  If G has a d-GUR framework then H has a d-GUR framework. \end{proposition}

\begin{proof} Let $G(p)$ be the $d$-GUR framework of $G$.  By Theorem 3.1 $G(p)$ has a PSD stress matrix $A$ of nullity $d+1$.  Let $H(p)$ be a $d$-dimensional framework of $H$ with the same vertex positions as $G(p)$.  $G(p)$ was generic so $H(p)$ is also generic.  Define the stress on the new added edge to be $0$ then $A$ is a stress matrix for $H(p)$ as well.  Thus $H(p)$ is a generic $d$-dimensional framework with a PSD stress matrix of nullity $d+1$.  Therefore $H$ has a $d$-GUR framework.  This proof is so trivial since adding an edge cannot do anything but make a rigid framework more rigid. \end{proof}  

\begin{theorem} Let G be a graph with $d+2$ or more vertices. Let H be obtained from G by a single $d$-Hennenberg operation.  If G has a d-GUR framework then H has a d-GUR framework. \end{theorem}

\begin{proof} Let $G(p)$ be the $d$-GUR framework of $G$.  By Theorem 3.1 $G(p)$ has a stress matrix $A$ which is PSD and of nullity $d+1$.  $G(p)$ must also be globally rigid, so by Theorem 2.4 it is redundantly rigid.  $G(p)$ is redundantly rigid, therefore it has a stress matrix $B$ for which the stress on every edge is non-zero. Pick $\epsilon_0$ so that every non-zero entry of $A$ is still non-zero in $A+\epsilon_0 B$. By Corollary 4.3, pick $\epsilon_1$ so that $A+\epsilon_1 B$ is still PSD and nullity $d+1$.  Pick $\epsilon=\Min(\epsilon_0, \epsilon_1)$ and let $\Omega_{G(p)}=A+\epsilon B$. $\Omega_{G(p)}$ is now an equilibrium stress matrix for $G(p)$ which is PSD of nullity $d+1$ and defines a non-zero stress $\omega_{ij}^p$ on every edge. Arrange $\Omega_{G(p)}$ so that the two vertices whose edge was eliminated by the $d$-Hennenberg operation, say $x$ and $y$, correspond to the last two columns and rows of the matrix respectively. Let $H(q)$ be a $d$-dimensional framework of $H$, where $q_i=p_i$ for all $i: \,  1\leq i \leq v$ and the extra vertex, $q_{v+1}=z$, of $H$ is placed somewhere along the line going through the eliminated edge $\{x,y\}$. Note that the configuration of $H(q)$ is almost the same as the configuration of $G(p)$, and that $H(q)$ is not a generic framework (see figure 2).  Let $\boldsymbol{xy}$ be the directed length of the old edge, let $\frac{1}{a}\boldsymbol{xy}$ be the directed length of edge $\{x,z\}$, and let $\frac{1}{b}\boldsymbol{xy}$ be the directed length of edge $\{z,y\}$ where we require $1/a+1/b=1$. Note choosing a value for $a$ or $b$ determines the placement of the final vertex $z$. We define a stress vector $\boldsymbol{\omega^q}$ on $H(q)$ as follows: 

(1) $\omega_{ij}^q=\omega_{ij}^p$ if $\{i,j\}$ is an edge of $G$ and $H$,

(2) $\omega_{xz}^q=a\cdot \omega_{xy}^p$, 

(3) $\omega_{zy}^q=b\cdot \omega_{xy}^p$, and

(4) $\omega_{zj}^q=0$ if $j \neq x,y$.  

These conditions constitute an equilibrium stress on $H(q)$. Let $\Omega_{H(q)}$ be the $ (v+1) \times (v+1)$ stress matrix defined by this stress. Arrange $\Omega_{H(q)}$ so that the first $v$ rows and columns correspond to the same vertices as in $\Omega_{G(p)}$, and the last row and column correspond to the new vertex $z$.  Let $\Omega_{G(p)}'$ be a $(v+1) \times (v+1)$ matrix, where the top left $v \times v$ block has the same entries as $\Omega_{G(p)}$ and the rest of the entries are zero. $\Omega_{G(p)}$ is PSD of nullity $d+1$, so $\Omega_{G(p)}'$ is PSD of nullity $d+2$. The last column and row of $\Omega_{G(p)}'$, which contains all zeros, corresponds to the vertex $z$.  Since the stress on each edge except for $\{x,y\}$, $\{x,z\}$, and $\{z,y\}$ is the same, $\Omega_{G(p)}'$ and $\Omega_{H(q)}$ have the same entries except for a $3 \times 3$ block in the bottom right corner.  Let $M$ be the $(v+1) \times (v+1)$ matrix such that $\Omega_{G(p)}' +M=\Omega_{H(q)}$. Let $M_{3 \times 3}$ be the non-zero $3 \times 3$ block of $M$.  Then
\begin{equation*}
M_{3 \times 3} = \omega_{xy} \left(
\begin{array}{ccc}
a-1 & 1 & -a \\
1 & b-1 & -b \\
-a & -b & a+b
\end{array} \right).
\end{equation*}      
Since $1/a+1/b=1$ it follows that $\Rank(M)=1$. If $\omega_{xy}>0$, pick $a,b=2$. Then the diagonal entries of $M$ are positive, so $M$ is PSD.  If $\omega_{xy}<0$, pick $a=-2, b=2/3$. Then the diagonal entries of $M$ are positive, so $M$ is PSD. The sum of two PSD matrices is PSD so $\Omega_{H(q)}$ is PSD.  Using Lemma 4.5, $\Ker(\Omega_{H(q)})=(\Ker(\Omega_{G(p)}') \cap \Ker(M))$.  Because $a, b \neq 0$, the standard basis vector $e_{v+1}$ is not in the kernel of $M$ but is clearly in the kernel of $\Omega_{G(p)}'$.  Thus, the dimension of the kernel of $\Omega_{H(q)}$ is strictly less than the dimension of the kernel of $\Omega_{G(p)}'$, so $d+1 \leq \Nullity(\Omega_{H(q)})< \Nullity(\Omega_{\emph{G(p)}}')=d+2$. Therefore, $\Omega_{H(q)}$ is PSD and of nullity $d+1$.  All that is left is to perturb $H(q)$ to a generic configuration.  To be able to do this $H(q)$ must be infinitesimally rigid, and the following proposition proves that it is.       

\begin{proposition}[Whiteley and Tay \cite{twir}] Suppose G(p) is a framework in $\RR^d$ and H is obtained from G by a Hennenberg operation. If H(q) is a framework such that the $d-1$ additional edges and the subdivided edge of G do not lie in a $d-1$ affine subspace, then H(q) is infinitesimally rigid. \end{proposition} 

In order for the $d-1$ additional edges and the subdivided edge to lie in a $(d-1)$-dimensional affine subspace $x$, $y$, $z$ and the $d-1$ vertices connected to $z$ must all lie in a $(d-1)$-dimensional affine subspace.  This means that for $G(p)$ the vertices $x$, $y$ and $d-1$ other vertices all lie in a $(d-1)$-dimensional affine subspace.  However such a configuration is not generic, contradicting our hypothesis on $G(p)$, so $H(q)$ is infinitesimally rigid.  

By Theorem 2.1, we now know that every generic framework of $H$ in dimension $d$ will be infinitesimally rigid, so the kernel of the transpose rigidity matrix for all of these frameworks will be minimal. Non-generic points are measure zero in $\RR^{dv}$ so $\forall \delta>0$ there exists a generic framework, $H(q')$ such that $|q'-q|<\delta$. The entries of the rigidity matrix are polynomial functions of the vertex positions, so by Lemma 4.4 $\forall \epsilon>0$ $\exists \delta>0$ such that the distance between $\Ker(df(H(q))^T)$ and $\Ker(df(H(q'))^T)$ is less than $\epsilon$. Let $\boldsymbol{v}=\boldsymbol{\omega^q}/\lVert \boldsymbol{\omega^q} \rVert$, then $\boldsymbol{v} \in \Ker(df(H(q))^T) \cap \mathbb{S}^{d-1}$.  Therefore $\exists \boldsymbol{u} \in \Ker(df(H(q'))^T) \cap \mathbb{S}^{d-1}$ such that $\lVert \boldsymbol{u}-\boldsymbol{v} \rVert < \epsilon$.   Define a stress $\boldsymbol{\omega^{q'}}$ on $H(q')$ such that $\boldsymbol{\omega^{q'}}=\boldsymbol{u}\lVert \boldsymbol{\omega^q} \rVert$ and let $\Omega_{H(q')}$ be the stress matrix defined by this stress vector. The stress on each edge of $H(q')$ differs from the stress on $H(q)$ by at most $\lVert \boldsymbol{\omega^q} \rVert\epsilon$.  Let $F$ be the matrix satisfying $F+\Omega_{H(q)}=\Omega_{H(q')}$.  The norm of the matrix $\lVert F \rVert_2=c\epsilon$ for some constant $c$. Finally, pick $q'$ so that $\epsilon$ is so small that F satisfies the hypothesis of Corollary 4.2.  Then $\Omega_{H(q')}$ is a PSD stress matrix of nullity $d+1$. By Corollary 3.4, $H(q')$ is a $d$-GUR framework. This completes the proof of the theorem. \end{proof}  

\begin{figure}[htb]
\begin{center}
\leavevmode
\includegraphics[width=0.8\textwidth]{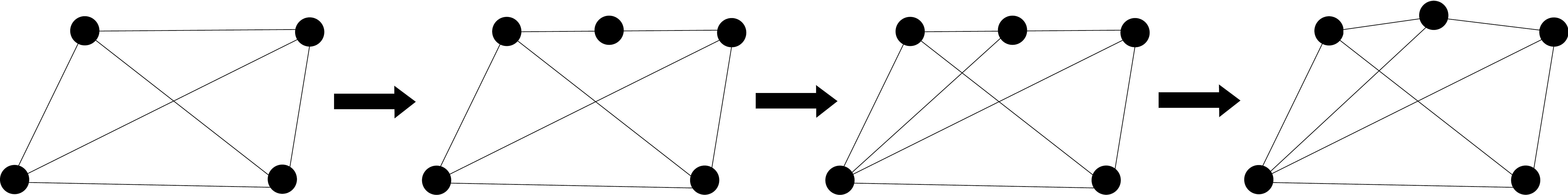}
\end{center}
\caption{Hennenberg operation on a edge with a positive stress.  The last step is a perturbation to a generic configuration.}
\label{figure2}
\end{figure}

\begin{figure}[htb]
\begin{center}
\leavevmode
\includegraphics[width=0.8\textwidth]{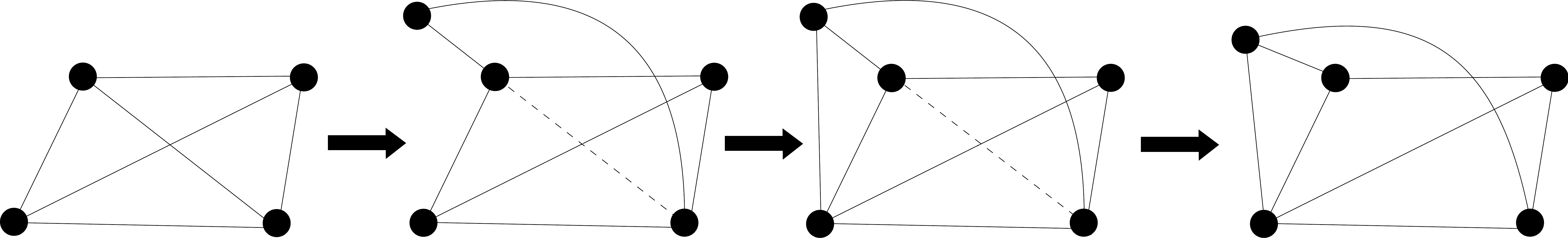}
\end{center}
\caption{Hennenberg operation on a edge with a negative stress.}
\label{figure3}
\end{figure}

\begin{corollary} Let G be a graph with $d+2$ or more vertices. Let H be obtained from G by any sequence of edge additions and $d$-Hennenberg operations. If G has a $d$-GUR framework then H does as well. \end{corollary}

\begin{corollary} For any d an infinite number of $d$-GGR graphs have a $d$-GUR framework.  \end{corollary}

\begin{proof} Any graph built from $d$-Hennenberg operations and edge additions to the complete graph $K_{d+2}$ will have a $d$-GUR framework by Corollary 5.4 and Proposition 2.5. \end{proof}

\textbf{Path, Path Addition, and Whitney Synthesis:} Let $G$ be a graph.  A \emph{path} in $G$ is a sequence of vertices of $G$ such that each vertex in the sequence is connected to the next vertex in the sequence by an edge.  A \emph{path addition} to a graph $G$ is the addition to $G$ of a path between two existing vertices of $G$, such that the edges and internal vertices of the path are not in $G$. A \emph{Whitney synthesis} of a graph $H$ from a graph $G$ is a sequence of graphs $G_0, G_1, \ldots, G_n$ where $G_0=G$ and $G_n=H$, and $G_i$ is the result of a path addition to $G_{i-1}$.

\begin{WST}[see e.g. Gross and Yellen, Chapter 5 Section 2 \cite{gy}] A graph G is 2-vertex connected if and only if G is a cycle or a Whitney synthesis from a cycle. \end{WST}

\begin{proposition}  Every 1-GGR graph with 3 or more vertices can be built from a sequence of Hennenberg operations and edge additions to $K_3$\footnote{Thanks to Timothy Sun for showing the Whitney Synthesis Theorem implies all 1-GGR graphs have a 1-GUR framework}. \end{proposition}  

\begin{proof} By Theorem 2.4 every $d$-GGR graph with $d+2$ or more vertices must be $d+1$-vertex connected and redundantly rigid.  For a 1-dimensional framework, redundant rigidity is equivalent to the graph being 2-edge-connected. 2-edge-connectedness is a weaker condition than 2-vertex-connectedness, so only 2-vertex-connectivity needs to be considered.  By the Whitney synthesis theorem every 2-vertex-connected graph is either a cycle or a path addition to a cycle. A 1-Hennenberg operation on a graph subdivides an edge by adding a vertex. A cycle with $v$ vertices can easily be generated from $K_3$ by performing $v-3$ Hennenberg operations on the graph. A path addition is simply an edge addition to a graph, followed by any number of Hennenberg operations on the added edge. Any sequence of path additions can be formed by successive edge additions and Hennenberg operations. \end{proof} 

\begin{theorem}[Connelly \cite{cn2}] If a graph $G$ with at least 4 vertices is redundantly rigid in $\RR^{2}$ and 3-vertex connected then $G$ can be obtained from $K_4$ by a sequence of Hennenberg operations and edge additions. \end{theorem} 

Redundant rigidity and 3-vertex connectivity are necessary and sufficient conditions for global rigidity in 2-dimensions, thus every 2-GGR graph on 4 or more vertices can be built from Hennenburg operations and edge additions on $K_4$

\begin{corollary} Every GGR graph in 1 and 2-dimensions has a universally rigid generic framework. \end{corollary}

\begin{proof} Suppose $G$ is a 2-GGR graph (1-GGR graph).  If $G$ has 4 (3) or more vertices then $G$ satisfies Theorem 5.7 (Proposition 5.6).  Thus Corollary 5.4 applies to $G$, so $G$ has a 2-GUR framework (1-GUR framework).  If $G$ has 3 (2) or fewer vertices than $G$ is a complete graph, and therefore universally rigid. This completes the proof of Theorem 1.1. \end{proof}

\section{Construction of an infinite family of $d$-SUR graphs} 

We now turn to the construction of an infinite family of $d$-SUR graphs. Recall that given some framework $G(p)$ the space of stresses on $G(p)$ is the kernel of the transpose rigidity matrix $df^T$. 

\begin{proposition} Let $G_0, G_1, \ldots, G_n$ be a sequence of graphs where $G_m$ is the result of a $d$-Hennenberg operation on $G_{m-1}$ and $G_0$ has a  $d$-GUR framework.  If $G_0(q)$ and $G_m(p)$ are generic frameworks then $\Nullity(df^T(G_m(p)))=\Nullity(df^T(G_0(q)))$.   \end{proposition}

\begin{proof}  By Theorem 5.2 and Proposition 2.6, each $G_m$ is $d$-GGR.  Therefore the generic frameworks of each $G_m$ are infinitesimally rigid so the rank of the rigidity matrix $df$ is maximal. Let $v_m$ and $e_m$ denote the number of vertices and edges of $G_m$ respectively. For a generic framework $G_m(p)$ $\Rank(df(G_m(p)))=dv_m-\binom{d+1}{2}$.  Let $C_m=\Nullity(df^T(G_m(p)))$ where $G_m(p)$ is generic. Then $C_m=e_m-dv_m+\binom{d+1}{2}$. Now we proceed by induction. The base case is trivial as $C_0=C_0$. Now suppose that $C_{m-1}=C_0$. Then \begin{equation}C_m=e_{m}-dv_{m}+\binom{d+1}{2}=d+e_{m-1}-(d+dv_{m-1})+\binom{d+1}{2}=C_{m-1}=C_0. \end{equation} Thus $C_m=C_0$.  \end{proof}

\begin{corollary} Let $G_0, G_1, \ldots, G_n$ be a sequence of graphs where $G_m$ is the result of a $d$-Hennenberg operation on $G_{m-1}$ and $G_0=K_{d+2}$ the complete graph on $d+2$ vertices.  If $G_m(p)$ is generic then $\Nullity(df^T(G_m(p)))=1$.  In particular $G_m(p)$ has a one dimensional space of stresses. \end{corollary}

\begin{proof} $G_0=K_{d+2}$ the complete graph on $d+2$ vertices thus $e_0=\binom{d+2}{2}$. Let $G_0(q)$ be a generic framework.  $G_0$ is $d$-GGR so $G_0(q)$ is infinitesimally rigid. Therefore \begin{equation} \Nullity(df^T(G_0(q))=e_0-dv_0+\binom{d+1}{2}=\binom{d+2}{2}-d(d+2)+\binom{d+1}{2}=1. \end{equation}  By Proposition 6.1 $\Nullity(df^T(G_m(p)))=1$ \end{proof}

\begin{theorem} Let $G_0, G_1, \ldots, G_n$ be a sequence of graphs where $G_0=K_{d+2}$ and $G_m$ is the result of a $d$-Hennenberg operation on $G_{m-1}$.  Then $\forall m>0$ the graph $G_m$ is $d$-SUR. \end{theorem}
\begin{proof} Each $G_m$ has a $d$-GUR framework. Thus we only need to show that for $m>0$ the graph $G_m$ has a generic $d$-dimensional framework which is not universally rigid. For some $l: 0 \leq l \leq n-1$ let $G_l(p)$ be a $d$-GUR framework of $G_l$.  By Corollary $G_l(p)$ has a $(d+l+2) \times (d+l+2)$ stress matrix $\Omega_{G_l(p)}$ which is unique up to scaling.  $G_l(p)$ is $d$-GUR so $\Omega_{G_l(p)}$ is PSD of nullity $d+1$.  $G_l(p)$ is also redundantly rigid so $\Omega_{G_l(p)}$ must also define a stress that is non-zero on every edge.  Let $G_m=G_{l+1}$ and set $k=(d+l+3)$.  We now proceed as in Theorem 5.2 and we obtain $k \times k$ stress matrices $\Omega_{G_l(p)}'$ and $\Omega_{G_m(q)}$.  The upper left $(k-1) \times (k-1)$ block of $\Omega_{G_l(p)}'$ has the same entries as $\Omega_{G_l(p)}$ and the entries of $k^{th}$ row and column of $\Omega_{G_l(p)}'$ are all zero.   $\Omega_{G_m(q)}$ and $\Omega_{G_l(p)}'$ have the same entries except for a $3 \times 3$ block in the bottom right corner.  $M$ is the $k \times k$ matrix such that $\Omega_{G_m(p)}' +M=\Omega_{G_m(q)}$. Let $M_{3 \times 3}$ be the non-zero $3 \times 3$ block of $M$.  Recall that
\begin{equation*}
M_{3 \times 3} = \omega_{xy} \left(
\begin{array}{ccc}
a-1 & 1 & -a \\
1 & b-1 & -b \\
-a & -b & a+b
\end{array} \right).
\end{equation*}      

Since $1/a+1/b=1$ it is not difficult to see that $\Rank(M)=1$. If $\omega_{xy}<0$ pick $a,b=2$. Then the diagonal entries of $M$ are negative so $M$ is NSD. If $\omega_{xy}>0$ pick $a=-2$ $b=2/3$. Then the diagonal entries of $M$ are negative so $M$ is NSD.  Let $\mathbf{e_{k}}$ be the last standard basis vector. \begin{equation} \mathbf{e_{k}}^T\Omega_{G_m(q)}\mathbf{e_{k}}=\mathbf{e_{k}}^T\Omega_{G_m(p)}'\mathbf{e_{k}} +\mathbf{e_{k}}^TM\mathbf{e_{k}}=0+ \omega_{xy}(a+b)=\omega_{xy}(a+b)<0. \end{equation}   Thus $\Omega_{G_m(q)}$ cannot be PSD.   However, $\Omega_{G_m(p)}'$ is a PSD $k \times k$ matrix of nullity $d+2$ so \begin{equation} \Rank(\Omega_{G_m(p)}')=k-(d+2)=l+1 \geq 1=\Rank(M). \end{equation} Thus it follows that $\Nullity(M) \geq \Nullity(\Omega_{G_m(p)}')$. The standard basis vector $e_k \notin \Ker(M)$ but $e_k \in \Ker(\Omega_{G_m(p)}')$ so the two kernels are not identical.  The dimension of the kernel of $M$ is greater than or equal to the dimension of the kernel of $\Omega_{G_m(p)}'$.  Therefore $\Ker(M)-\Ker(\Omega_{G_m(p)}') \neq \varnothing$.  So $\exists \boldsymbol{v} \in \Ker(M)$ such that $\boldsymbol{v} \notin \Ker(\Omega_{G_m(p)}')$.   \begin{equation} \boldsymbol{v}^T\Omega_{G_m(q)}\boldsymbol{v}=\boldsymbol{v}^T\Omega_{G_m(p)}'\boldsymbol{v} +\boldsymbol{v}^TM\boldsymbol{v}=\boldsymbol{v}^T\Omega_{G_m(p)}'\boldsymbol{v}+0=\boldsymbol{v}^T\Omega_{G_m(p)}'\boldsymbol{v}>0. \end{equation}  Thus $\Omega_{G_m(q)}$ cannot be NSD. Therefore $\Omega_{G_m(q)}$ is indefinite.  Again we proceed as in Theorem 5.2 and we obtain a generic framework $q'$ such that $\Omega_{G_m(q')}-\Omega_{G_m(q)}=F$ and $\lVert F \rVert=c\epsilon$.  We then pick $\epsilon$ small enough that F satisfies the hypothesis of Corollary 4.1.  Then $\Omega_{G_m(q')}$ has at least as many positive and negative eigenvalues as $\Omega_{G_m(q)}$.  Since $\Omega_{G_m(q)}$ is indefinite $\Omega_{G_m(q')}$ must also be indefinite. By Corollary 6.2 $\Omega_{G_m(q')}$ is a unique stress matrix for $G_m(q')$ up to scale.  Therefore the stress matrices of $G_m(q')$ are either indefinite or are the zero matrix.  The zero matrix has nullity greater than $d+1$ so $G_m(q')$ has no PSD stress matrix of nullity $d+1$.  Finally by Corollary 3.4 $G_m(q')$ cannot be universally rigid.  So $G_{m}$ is a $d$-SUR graph. \end{proof}

\begin{corollary} Let G be a 2-GGR graph with more than 4 vertices.  If the generic frameworks of G have a one dimensional space of stresses then G is 2-SUR. \end{corollary}  
\begin{proof} The statement is a direct result of Theorem 5.7 and Theorem 6.3. \end{proof}

\end{document}